%
%

\magnification=1200

\font\titfont=cmr10 at 12 pt

\font\headfont=cmr10 at 12 pt



\def\AAA{\bf}

\overfullrule=0in

\def\boxit#1{\hbox{\vrule
 \vtop{%
  \vbox{\hrule\kern 2pt %
     \hbox{\kern 2pt #1\kern 2pt}}%
   \kern 2pt \hrule }%
  \vrule}}

  \def\harr#1#2{\ \smash{\mathop{\hbox to .3in{\rightarrowfill}}\limits^{\scriptstyle#1}_{\scriptstyle#2}}\ }

 \def\GG{{{\bf G} \!\!\!\! {\rm l}}\ }

\def\bra#1#2{\langle #1, #2\rangle}

\def\ss{\subset}

\def\half{\hbox{${1\over 2}$}}

\def\max{{\rm max}}

\def\Sym{{\rm Sym}^2}

\def\rn{\bbr^n}

\def\Int{{\rm Int}}

\def\Symn{{\Sym(\rn)}}

\def\Theorem#1{\medskip\noindent {\bf THEOREM \bf #1.}}

\def\Lemma#1{\medskip\noindent {\bf Lemma #1.}}
\def\Remark#1{\medskip\noindent {\bf Remark #1.}}
\def\Note#1{\medskip\noindent {\bf Note #1.}}
\def\Def#1{\medskip\noindent {\bf Definition #1.}}

\def\Ex#1{\medskip\noindent {\bf Example \bf    #1.}}

\def\pf{\medskip\noindent {\bf Proof.}\ }
\def\qed{\hfill  $\vrule width5pt height5pt depth0pt$}

   \def\cp{{\cal P}}

\def\cp{{\cal P}}

\def\wt{\widetilde}

\def\and{\qquad {\rm and} \qquad}

\def\bbr{{\bf R}}

\def\d{\delta}
\def\e{\epsilon}

\def\l{\lambda}

\def\O{\Omega}

\def\bo{\partial \Omega}

\def\Symn{\Sym(\rn)}
 
\def\USC{{\rm USC}}

\def\cpt{\wt{\cp}}
\def\ft{\wt F}
\def\ob{\overline{\O}}

\def\AAA{1} \def\AA{2}
\def\BB{3} \def\CC{4}  \def\EE{5}\def\HH{6} \def\FF{7}\def\GG{8}

\centerline
{
\titfont THE AE THEOREM AND  ADDITION THEOREMS }

\smallskip

\centerline
{
\titfont FOR QUASI-CONVEX FUNCTIONS }
\smallskip

\bigskip

\centerline{\titfont F. Reese Harvey and H. Blaine Lawson, Jr.$^*$}
\vglue .9cm
\smallbreak\footnote{}{ $ {} \sp{ *}{\rm Partially}$  supported by
the N.S.F. } 

\vskip .2in

\centerline{\bf ABSTRACT} \medskip
  \font\abstractfont=cmr10 at 10 pt

  {{\parindent= 1.0in

\narrower\abstractfont \noindent
The main point of this paper is to prove the following useful result:
If the almost everywhere 2-jet of a locally quasi-convex function $u$
satisfies a degenerate elliptic constraint $F$, then $u$ is $F$-subharmonic,
i.e., $u$ is a viscosity $F$-subsolution. This AE Theorem makes 
otherwise difficult results transparent. Some instances of this
 are presented, including two versions of addition, and 
a comparison theorem.

}}

\vskip.5in


\centerline{\headfont \AAA.\ Introduction.}
 \medskip 

The main point of this short article is to state and prove a highly useful theorem
giving a necessary and sufficient condition for a quasi-convex function to be
a subsolution of a degenerate elliptic equation. If the equation is  pure second-order 
with constant coefficients, the result was established in [HL$_1$, Cor.7.5],
and here it is generalized to the fullest extent. As in out previous work
(cf. [HL$_{1,2}$]), we focus on viscosity subsolutions (rather than solutions)
and take a ``potential theoretic approach'' in parallel with, and in fact generalizing,
the theory of plurisubharmonic functions in several complex variables.
More specifically, fix $X^{\rm open}\ss \rn$, or more generally a manifold, and consider a closed subset
$$
F \ \ss\ J^2(X) \equiv X\times\bbr\times \rn\times\Symn
$$
of the bundle of 2-jets of functions on $X$ which satisfies the mild
positivity condition (or degenerate ellipticity condition) (P) below.
Using  test functions, one defines the notion of {\sl $F$-subharmonicity}
(or {\sl viscosity $F$-subsolution}) for any upper semi-continuous function
$u:X\to [-\infty,\infty)$. The set $F(X)$ of these functions enjoys a majority
of the important properties of the family of classical subharmonic functions
on $X$.

Recall now that a quasi-convex function on $X$ is twice differentiable at
almost every point.  The main result is:

\medskip
\noindent
{\bf The AE Theorem.}  {\sl
Suppose $u$  is locally quasi-convex on a manifold $X$.    Then
\medskip
\centerline
{
$J^2_x u \in F$ \ \ a.e.   $ 
\qquad \Rightarrow\qquad
u\in F(X)$
}
}
\medskip

This theorem has many useful applications, and while
 perhaps not a surprise to experts, it deserves to be highlighted.
 (We were unable to find it in the literature outside the special case in [HL$_1$].)
It proved, for example, to be a convenient tool in [HL$_3$]
and in characterizing  radial subharmonics in [HL$_5$], a fact which prompted the authors to write
this note.  

The result applies immediately to prove {\sl addition}. Suppose $F,G \ss J^2(X)$
are closed subsets satisfying (P) as  above and let $H\equiv \overline{F+G}$ (fibre-wise sum).

\medskip
\noindent
{\bf The Basic Addition Theorem.}  {\sl Suppose $u$ and $v$ are locally quasi-convex on
$X$.  Then
\medskip
\centerline
{
$u\in F(X) \ \ {\rm and}\ \ v\in G(X)
\qquad \Rightarrow\qquad
u+v\in H(X)$
}
\noindent
Moreover, if $F$ and $G$ are constant coefficients subequations in $\rn$, the quasi-convex assumptions can be dropped.}
\medskip

A nice application in the constant coefficient case is provided by the comparison result 
Theorem \HH.4 and its generalization Theorem \HH.8.
A stronger form of the addition theorem is proved in \S \FF.

\vskip.3in


\centerline{\headfont \AA.\ Preliminaries.}
 \medskip 

There are several (equivalent) ways of defining subsolutions.  For the purposes
of this paper the most convenient approach is as follows.   Let $\Symn$ denote the
set of symmetric $n\times n$ matrices.

\Def{\AA.1}  Given a real-valued function $w$ defined on an open subset $X\ss\rn$,
a point $x\in X$ is an {\bf upper contact point for $w$}
if there exists $(p,A) \in\rn\times\Symn$ such that
$$
w(y)\ \leq\ w(x) + \bra p {y-x} +\half \bra{A(y-x)}{y-x} \qquad\forall\, y \ {\rm near}\  x.
\eqno{(\AA.1)}
$$
In this case, $(p,A)$ is called an {\bf upper contact jet for $w$ at $x$}.

Given a subset  $F\ss X\times \bbr\times \rn\times \Symn \equiv J^2(X)$,
we adopt the following definition of {\sl subsolution}.
Let $X\ss\rn$ be an open subset, and denote by $\USC(X)$ the set of 
upper semi-continuous functions taking values in $[-\infty, \infty)= \bbr\cup \{-\infty\}$.
Let $F_x$ denote the fibre of $F$ at $x$.

\Def{\AA.2} A function $w\in \USC(X)$ is said to be $F$-{\bf subharmonic on $X$}
if for each point $x\in X$ and every upper contact jet $(p,A)$ for $w$ at $x$,  the jet
$(w(x),p,A) \in F_x$.   Let $F(X)$ denote the set of functions
which are $F$-subharmonic on $X$.
\medskip

Conditions must be placed on $F$ in order for this concept to be useful. 
Of crucial importance is the condition of {\bf positivity}:
$$
(x,r,p,A) \in F 
\qquad\Rightarrow\qquad
(x,r,p,A+P)\in F \ \ \ \forall \ P\geq0.
\eqno{(P)}
$$
It  is essential in order  for every $C^2$-function $w$
with $(x,w(x), D_xw, D^2_x w)\in F$ for all $x\in X$ to
be subharmonic. The second condition:
$$
F \ \ {\rm is\ a\ closed\ subset}
\eqno{(C)}
$$
is a necessary requirement for the elementary properties
of $F(X)$ involving limits and upper envelopes  to hold
(see Theorem 2.6 in [HL$_2$]).  

If $F$ satisfies these minimal conditions (P) and (C)
we say that $F$ is a  {\bf primitive subequation}.
The results of this paper hold under (P) and (C) with the exception
of comparison results where $F$ is required to be a subequation (see [HL$_2$]).

There are two extreme cases where the set of upper contact jets of a function $w$
at a point $x$ are 
essentially completely understood. The first case is where $w$ has no upper contact
jets at $x$. Such is the case if $w(x)=-\infty$ for example.  It is also true of the function
$w(x)=|x|$ at $x=0$.  The other extreme is when $w$ is twice differentiable at $x$.
By basic differential calculus we have the following.

\Lemma{\AA.3} 
{\sl Suppose $w$ is twice differentiable at $x$.  If $(p,A)$ is an upper contact
jet for $w$ at $x$, then $p=D_xw$ is unique and $A=D^2_xw +P$ for some $P\geq0$.
Conversely, for each $ P>0$ $(D_xw, D^2_xw + P)$ is an upper contact jet for $w$ at $x$.
}

\medskip

Adding a smooth function $\psi$ to $w$ does not change the set of upper contact
points but only the upper contact jets.

\Lemma{\AA.4}  {\sl
Suppose $\psi$ is smooth.  Then
$$
x\ \ {\sl is\  an\  u.\ c. \ point\ for\ \ } w
\quad\iff\quad
x\ \ {\sl is \ an\ u. \  c. \ point\ for\ \ } w+\psi
$$
$$
(p,A) \ \ {\sl is\  an\ u. \  c. \ jet\ for\  } w \ {\sl at\ } x 
\quad\iff\quad
(p+ D_x\psi, A+ D^2_x\psi) \ \ {\sl is\  an\ u. \ c. \  jet\ for\  } w+\psi \ {\sl at\ } x
$$
}

A third elementary fact needed here concerns convex functions.

\Lemma{\AA.5} {\sl
If $w$ is convex and twice differentiable  at $x$, then $D^2_xw\geq0$.
}

\Ex{\AA.6}  Each function $w\in \USC(X)$ determines a smallest primitive subequation,
denoted $\overline{J^+(w)}$,  with the property that  $w$ is subharmonic.  Namely, let $J^+(w)$
denote the set of tuples $(x, w(x), p,A)$ such that $(p,A)$ is an upper contact jet for $w$ at $x$, 
and then take the closure $\overline{J^+(w)}$ in $X\times \bbr\times\rn\times\Symn$.
  Obviously $J^+(w)$ satisfies condition (P), and 
hence the closure also satisfies (P).  This example $\overline{J^+(w)}$ is in some sense 
pathological since it does not satisfies the requirements of negativity or the topological
regularity required of a subequation in [HL$_2$].

\medskip

Recall that a function $w$ is called   $\l$-{\bf quasi-convex} if the function
$$
u(y) \ \equiv \ w(y) + {\l \over 2} |y|^2
\eqno{(\AA.2)}
$$
is convex, and $w$ is called {\bf quasi-convex} if this is true for some $\l\geq0$.
 Note that the concept of 
being locally quasi-convex is preserved under diffeomorphisms
and therefore makes sense on any manifold.

Any property of convex functions can always be reformulated via (\AA.2)
as a property of quasi-convex functions. Conversely, each result concerning
quasi-convex functions  can  be reformulated as a result about convex functions.

The possible upper contact jets of a (locally) quasi-convex function are fairly well understood.
We discuss this  in the next section.

\vskip .3in


\centerline{\headfont \BB.\ Upper Contact Jets of Quasi-Convex Functions}
 \bigskip 

Elementary proofs of the  results in this section can be found, for example, in the notes
[HL$_4$] posted on the ArXiv.\medskip

\centerline{\bf First Derivatives}\medskip

For quasi-convex functions, differentiability is automatic at each upper contact point.

\Lemma{\BB.1. (D at UCP)} {\sl Suppose $w$ is quasi-convex.  If $x_0$ is an upper contact point
for $w$, then $w$ is differentiable at $x_0$. Moreover, if $(p,A)$ is any upper contact
jet for $w$ at $x_0$, then $p=D_{x_0}w$ is unique.
}
\medskip

Another even more standard result is called partial continuity of the gradient, or first derivative.

\Lemma{\BB.2. (PC of FD)} {\sl Suppose $w$ is quasi-convex and $x_j \to x_0$.
If $w$ is differentiable at each $x_j$ and at $x_0$, then $D_{x_j}w\to D_{x_0}w$.
}

\bigskip

\centerline{\bf Second Derivatives}\medskip

The results concerning the second-order part of upper contact jets of quasi-convex functions
are of a deeper nature.  The almost everywhere existence of the first derivative 
was omitted from the previous discussion because of the following stronger result.

\Theorem {\BB.3. (Alexandrov)} {\sl
 A locally quasi-convex function is twice differentiable almost everywhere.
}
\medskip

For the next result we need two variations of the notion of an upper contact jet.
First, we say that $(p,A)$ is a {\bf strict} upper contact jet for $w\in \USC(X)$ at $x_0\in X$ if
the upper contact inequality (\AA.1) is strict for $y\neq x_0$.
An understanding of the strict  upper contact jets will be  adequate for 
our discussion since $(p,A)$ is an upper contact jet 
 if and only if   $(p, A+\e I)$ is a strict upper contact jet for all $\e>0$.
Second, we need a   notion of upper contact point and jet, which requires the inequality
(\AA.1) to hold globally.

\Def{\BB.4}  Given $w\in \USC(X)$ and $A\in \Symn$, a point $x$ is called a
{\bf global upper contact point of type $A$ on $X$} if for some $q\in\rn$
$$
w(y)\ \leq\ w(x) + \bra q {y-x}  + \half\bra{A(y-x)}{y-x}\quad \forall\, y\in X.
\eqno{(\BB.1)}
$$
Let $C(w,X,A)$ denote the set of all global  upper contact points of type $A$ on $X$ for the function $w$.

\Remark{\BB.5} Note that if $w$ is quasi-convex, then by (D at UCP) each point
$x\in C(w,X,A)$ is a point of differentiability and the only $q$ in (\BB.1) is $q=D_xw$.

\Theorem{\BB.6. (Jensen-Slodkowski)} {\sl
Suppose that $w$ is a quasi-convex function possessing a strict upper contact
jet $(p_0,A_0)$ at $x_0$.  
Let $B_\rho$ denote the ball of radius $\rho$ about $x_0$.
Then there exists $\bar \rho >0$ such that the measure}
$$
\left|   C(w, B_\rho, A_0)      \right| \ >\ 0 \qquad \forall\, 0 < \rho\leq \bar \rho.
\eqno{(\BB.2)}
$$

The four results above yield the following  theorem which is 
easy to apply and adequate for many purposes.  It is concerned with the
upper contact jets of a quasi-convex function.
The order two part of this theorem can 
be considered a ``partial upper semi-continuity of the
second derivative''.  We will use the acronym PUSC of SD.

\Theorem{\BB.7. (Upper Contact Jets)}
{\sl 
Suppose $w$ is quasi-convex with an upper contact jet $(p_0,A_0)$ at a point $x_0$.
 Then
\medskip
{\rm (D at UCP)}  \quad  $w$ is differentiable at $x_0$ and $D_{x_0}w = p$.
\medskip

Suppose $E$ is a set of full measure in a neighborhood of $x_0$. Then
there exists a sequence $\{x_j\} \ss E$ with $x_j\to x_0$  such that $w$ is twice differentiable
at each $x_j$ and 
\medskip
{\rm (PC of FD)}  \quad $D_{x_j} w \ \to\ D_{x_0} w=p_0$,
\medskip
 
{(\rm PUSC of SD)}  \quad $D_{x_j}^2 w \ \to\ A \leq A_0$.
}

\pf
By Alexandrov's Theorem, the set of points $x\in E$ where $w$ is twice differentiable,
is a set of full measure. In order to apply the Jensen-Slodkowski Lemma we replace 
$(p_0, A_0)$ by the strict upper contact jet $(p_0, A_0+\e I)$. Now choose a sequence
$\e_j\to0$, and pick a point $x_j \in B_{\e_j}(x_0)$ such that: (1) $x_j\in E$, (2) $w$ is twice differentiable at $x_j$, and (3) $x_j$ is a global upper contact point of type $A_0+\e_j I$
on $B_{\e_j}(x_0)$ for $w$.  By the basic differential calculus fact 
Lemma \AA.3, $D_{x_j}^2w \leq A_0+\e_j I$. 
Since $w$ is $\l$-quasi-convex, we have $D^2_{x_j} w +\l I \geq0$. Thus,
$$
-\l I \ \leq \ D^2_{x_j} w\ \leq\ A_0+\e_j I.
\eqno{(\BB.3)}
$$
By compactness there is  a subsequence such that $D^2_{x_j} w \to A\leq A_0$.\qed
\medskip

Theorem \BB.7 can be stated succinctly in terms of the subset $J^+(w)$ (see Example \AA.6)
representing the upper contact jets of $w$, and another  subset depending on $E$
Define the subset $J(w,E)$ of the 2-jet bundle
$J^2(X)$ to be the set of tuples $(x,w(x), D_x w, D^2_x w+P)$ such that $x\in E$, $w$ 
is twice differentiable at $x$, and $P\geq0$.  Then Theorem \BB.7  condenses to:
$$
{\rm If\ } w \ {\rm is \ quasi\!-\!convex\ and \ } E\ {\rm has\ full\  measure,\ then\ } \ J^+(w) \ \ss\ \overline{J(w,E)}.
\eqno{(\BB.4)}
$$

\vfill\eject


\centerline{\headfont \CC.\ The Almost Everywhere Theorem.}
 \medskip

 \medskip
The main result of this section is the following.
If $u$ is twice differentiable at $x$, let $J^2u$ denote the 2-jet $(u(x), D_xu, D^2_xu)$ at $x$.

\Theorem{\CC.1} {\sl
Suppose that $F$ is a primitive subequation on a manifold $X$ and $u:X\to\bbr$
is a locally quasi-convex function. Then}
$$
J_x^2 u \in F_x \ {\sl for\ almost\ all\  } x\in X \qquad \Rightarrow
\qquad u \ {\sl is \ } F\,{\sl subharmonic\  on \ } X.
\eqno{(\CC.1)}
$$
\pf
If $u$ is not $F$-subharmonic on $X$, then there exists an upper contact jet $(p_0,A_0)$ for 
$u$ at some point $x_0\in X$ with $(u(x_0),p,A) \notin F_{x_0}$.

We now apply Theorem \BB.4 to $u$ with $E$ taken to be the set of second
differentiability points $x$ for $u$ where the  2-jet of $u$ belongs to $F_x$,
i.e., $J^2_x u \in F_x$.
This yields a sequence $x_j\to x_0$ with $J^2_{x_j} u  \to (u(x_0), p_0, A_0-P)$ with 
$P\geq0$.  The fact that $F$ is closed, along with the positivity condition, implies that
$(u(x_0),p_0,A_0) \in F_{x_0}$, which is a contradiction.\qed
\medskip

By  the elementary calculus Lemma \AA.3 and Alexandrov's Theorem
 the concluding implication  (\CC.1) can be replaced by the ``if and only if'' statement:
$$
J^2_x u \in F_x \ \ {\rm a.e.}
\qquad\iff\qquad
u\in F(X)
\eqno{(\CC.1)'}
$$

\vskip .3in


\centerline{\headfont \EE.\  The Basic Addition Theorem}
 \medskip 

\Theorem{\EE.1} {\sl
Suppose that $F$  and $G$ are primitive subequations on a manifold $X$. Then
$H\equiv \overline{F+G}$ is also a primitive subequation.
Moreover, if  $u$ and $v$ are quasi-convex functions on $X$, then}
$$
u\in F(X)\ \ {\sl and}\ \ v\in G(X)
\qquad\Rightarrow\qquad
u+v\in H(X).
$$
Moreover,  If $F$ and $G$ are constant coefficient (translation invariant)
 on an open subset $X\ss \rn$,
then the assumption that $u$ and $v$ are quasi-convex can be dropped.
\pf 
First note that the sum $F+G$ satisfies condition (P) and therefore so does its closure.
Now by Alexandrov's Theorem and Lemma \AA.3,
 $J^2_x u \in F_x$ and $J^2_x v \in G_x$ for  a.a. $x\in X$.
  Hence $J^2_x (u+v) = J^2_xu+J^2_xv$ belongs to  $H_x$ for a.a. $x$. Now the AE Theorem \CC.1 
  applied to $H$ shows that  $u+v\in H(X)$.

The constant coefficient part  follows easily  from quasi-convex approximation (cf. [HL$_1$]).
\qed

\vfill\eject


\centerline{\headfont \HH.\  Comparison in the Constant Coefficient Case}
 \medskip 
 
 \centerline{ \bf Case I -- Pure Second-Order}
 \medskip
 
 In this case comparison always holds, and this follows in a straightforward 
 intuitive manner from the Addition Theorem \EE.1 along with some
 algebraic considerations.
 
 We note that a {\sl constant coefficient pure second-order subequation}
 is any closed subset $F\ss\Symn$ which satisfies
 $$
 A\in F \qquad\Rightarrow\qquad A+P \in F\quad \forall\, P\geq0.
 \eqno{(P)}
 $$
 The most basic of all examples is
 $$
 \cp\ :\ \ \l_{\rm min}(A)\ \geq\ 0
 $$
 where $\l_{\rm min}(A)$ is the minimum eigenvlaue.  The dual subequation is
 $$
 \cpt\ :\ \ \l_{\rm max}(A)\ \geq\ 0
 $$
 the so-called {\sl subaffine (or co-convex) subequation}.   The name subaffine for the subequation 
 $\cpt$ is justified by the following discussion.
 The reader is invited to verify the following facts (or see the proofs of Lemma 3.3a and Corollary 3.5 in 
 [HL$_2$]):
 $$
 \quad\quad
 F+\cp \ \ss\ F 
 \qquad\iff\qquad
 \ft +\cp \ \ss\ \ft \qquad{\rm and}
 \eqno{(\HH.1)}
 $$
 $$
 F+\ft \ \ss\ \cpt.
 \eqno{(\HH.2)}
 $$
  Let  Aff denote the space of all affine functions on $\rn$.
  A function $w\in\USC(X)$ will be called a {\sl subaffine function on $X$} if
 $$
 w \leq a \ \ \ {\rm on}\ \ \partial K 
  \qquad\Rightarrow\qquad
   w \leq a \ \ \ {\rm on}\ \ K \qquad \forall\, K^{\rm cpt} \ss X \ \ {\rm and}\ \ \forall \,a\in {\rm Aff}
 \eqno{(\HH.3)}
 $$
 
The final step in the proof of comparison requires the following characterization of $\cpt$-subharmonic functions.
  
  \Lemma{\HH.1. ([HL$_1$])}
  $$
  w\in \cpt(X) 
  \qquad\iff\qquad
  w \ \ {\sl is\  a\ subaffine\ function \ on\ \ } X.
  $$

  Now putting these results together we have two theorems.
  
  \Theorem{\HH.2. (The Subaffine Theorem)}
  {\sl
  If $u\in F(X)$ and $v\in\ft(X)$, then }
  $$
{\sl the\ sum\ \ } w= u+v \ \ {\sl is\ a\ subaffine\  function\ on\ \ } X.
  $$
  \pf By  Theorem \EE.1 and (\HH.2), $w$ is $\cpt$-subharmonic so that Lemma \HH.1 applies.\qed
  
  \medskip
  
  In order to state comparison for an arbitrary domain $\O\ss\rn$, 
  (i.e., a bounded connected open subset), we note the following 
  strengthening of  (\HH.3).  To give the reader some motivation, we point to the examples mentioned in 
  Remark \HH.9, which arise in studying the Dirichlet Problem with 
  Prescribed Singularities [HL$_6$].  We note, however, that this
  lemma is required even when $\O$ is the open ball.
  
  \Lemma{\HH.3} 
  {\sl
  Let $\O\ss\rn$ be any domain, and set
  $\bo \equiv \ob-\O$.  Then
  $$
    w \leq a \ \ \ {\rm on}\ \ \partial \O 
  \qquad\Rightarrow\qquad
   w \leq a \ \ \ {\rm on}\ \ \ob
   \eqno{(\HH.4)}
 $$
 for all functions $w\in\USC(\ob)$ which are subaffine on $\O$.
  }
  \pf
  We can assume $a\equiv 0$. Exhaust $\O$ by compact sets $K_1 \ss K_2\ss \cdots$ and 
  set $U_\d \equiv \{x\in\ob : w(x) < \sup_{\bo} w + \d\}$.
  Since $w \in \USC(\ob)$, each $U_\d$ is an open neighborhood of $\bo$ in $\ob$.
  Therefore, for $j$ sufficiently large we have $\partial K_j \ss U_\d$. Hence by (\HH.3),
  for all  $j$ sufficiently large we have
  $$
  \sup_{K_j} w \ \leq\   \sup_{U_\d} w  \ \leq\   \sup_{\bo} w +\d
  $$
  proving that $\sup_{\O} w\leq \sup_{\bo} w +\d$ for all $\d>0$.\qed

 \Theorem{\HH.4. (Comparison)}
  {\sl
Suppose $\O\ss\rn$ is a domain, and $u,v \in \USC(\ob)$.  
If $u\in F(\O)$ and $v\in \ft(\O)$, then}
$$
u+v\ \leq 0 \quad{\rm on}\ \ \bo
\qquad\Rightarrow\qquad
u+v\ \leq\ 0 \quad{\rm on}\ \ \ob.
$$
  
  \pf
  Apply Theorem \HH.2 and Lemma \HH.3, using the affine function $a\equiv 0$ in (\HH.4).\qed
  \medskip
  
  Both of these two  theorems were first proved in [HL$_1$].

  For simplicity and clarity we discussed Case I separately,
  although it is a subset of the next case.
  
  \vskip.3in
  
 \centerline{ \bf Case II -- Gradient-Free Subequations}
 \medskip
 
 Again in this case comparison always holds.
 
 \Def{\HH.5} A {\sl gradient-free constant coefficient subequation} is a closed subset
 of $\bbr\times \Symn$ satisfying {\sl positivity} and {\sl negativity}:
 $$
 (r,A) \in F \qquad\Rightarrow\qquad (r-s, A+P)\in F\qquad \forall\, s\geq0\ \ {\rm and}\ \ \forall\, P\geq0.
 \eqno{(P)\  {\rm and} \  (N)}
 $$
 
 Now the most basic example  is the subequation
 $$
 \bbr_-\times \cp\qquad{\rm where}\ \ \bbr_-\ \equiv\ (-\infty, 0].
 $$
  The dual subequation is
  $$
  \wt {\bbr_-\times \cp} \ =\ (\bbr_-\times \Symn) \cup (\bbr\times \cpt),
 \eqno{(\HH.5)}
 $$
  that is,
 $$
  \wt {\bbr_-\times \cp}  \ :\   \  r\leq 0 \ \ \ {\rm or}\ \ \ \l_{\rm max}(A)\ \geq\ 0.
 \eqno{(\HH.5)'}
 $$
 
 Here the important algebraic facts are as follows. For any set $F\ss \bbr\times \Symn$:
 $$
F\ \ {\rm is\ a\ subequation} 
 \qquad\iff\qquad
\ft\ \ {\rm is\ a\ subequation} 
  \eqno{(\HH.6)}
 $$
in which case
 $$
 F + \ft \ \ss\   \wt {\bbr_-\times \cp} 
  \eqno{(\HH.7)}
 $$
 Let Aff$^+$ denote the space of functions of the form
 $$
 a^+\ =\ \max\{a,0\}\qquad{\rm where} \ \ a \in  {\rm Aff}\ \ {\rm is\ affine}.
 $$
 These will be referred to as the {\sl affine plus functions} on $\rn$.
 The enhancement of Lemma \HH.1 characterizes the $(\wt {\bbr_-\times \cp})$-subharmonics.
 
 \Lemma{\HH.6}
 {\sl
One has  $w\in (  \wt {\bbr_-\times \cp} )(X) \quad\iff\quad w\in\USC(X)$ and}
 $$
 w \leq a^+ \ \ \ {\rm on}\ \ \partial K
  \qquad\Rightarrow\qquad
   w \leq a^+ \ \ \ {\rm on}\ \ K \qquad \forall\,  K^{\rm cpt} \ss X \ \ {\rm and} \ \ \forall\, a^+\in {\rm Aff}^+ 
 $$
 
Such functions will be referred to as {\sl subaffine plus functions on $X$}.
 Here as before, if $\O\ss \rn$ is any domain and $w$ is any upper semi-continuous function on $\ob$ 
 which is subaffine plus on $\O$, then
 $$
 w\ \leq \ a^+ \ \ \ {\rm on}\ \ \bo
 \qquad\Rightarrow
 \qquad
 w\ \leq \ a^+ \ \ \ {\rm on}\ \ \ob \quad \forall\, a^+\in {\rm Aff}^+.
 \eqno{(\HH.8)}
 $$

 The enhancements of Theorems \HH.2 and \HH.4 are now immediate.
 We assume $F\ss \bbr\times\Symn$ is a gradient-free  subequation.
 
 \Theorem{\HH.7 (The Subaffine Plus Theorem)}
{\sl
 If $u\in F(X)$ and $v\in \ft(X)$, then}
  $$
  {\sl the\ sum\ }\ w\equiv u+v \ \ {\sl is\ subaffine\ plus\ on\ \ } X.
  $$
  
  \pf Apply  the Addition Theorem \EE.1,  (\HH.7) and Lemma \HH.6.   \qed
  
 \Theorem{\HH.8. (Gradient-Free Comparison)}
  {\sl
Suppose $\O\ss\rn$ is a domain, and $u,v \in \USC(\ob)$.  
If $u\in F(\O)$ and $v\in \ft(\O)$, then}
$$
u+v\ \leq 0 \quad{\rm on}\ \ \bo
\qquad\Rightarrow\qquad
u+v\ \leq\ 0 \quad{\rm on}\ \ \ob.
$$

 \pf Apply Theorem \HH.7 and then (\HH.8).  Note that $0\in$ Aff$^+$ is an affine plus function.\qed
 
 \Remark {\HH.9}  If $\O$ is a domain and $x_1,...,x_k\in \O$, then $\O-\{x_1,...,x_k\}$ is also a 
 domain, and the theorems above apply to this domain as well. This is particularly useful in our
 study of the Dirichlet Problem with Prescribed Singularities [HL$_6$] on a domain $\O$ with smooth boundary.

\vfill\eject


\centerline{\headfont \FF.\  A Stronger Form of the  Addition Theorem}
 \medskip

\Theorem{\FF.1}
{\sl
Suppose $u$ and $v$ are quasi-convex and that  the sum $w\equiv u+v$
has an upper contact jet $(p_0,A_0)$ at $x_0\in \rn$.
Then $x_0$ is an upper contact point for both $u$ and $v$.  Furthermore:
\medskip
\centerline{\rm First Derivatives}
\medskip
\noindent
{\bf (D at UCP)}.  Both $u$ and $v$ are differentiable at $x_0$, 
and the upper contact jets are all of the form $(D_{x_0}u, -)$ and 
$(D_{x_0}w, -)$ respectively.
\medskip
\noindent
{\bf (DPC at FD)}. For any sequence $x_j\to x_0$ with both $u$ and $v$ 
differentiable at each $x_j$,
\medskip
\centerline{$D_{x_j}u \to D_{x_0}u
\and
D_{x_j}v \to D_{x_0}v$}
\medskip
\centerline{\rm Second Derivatives}
\medskip
\noindent
{\bf (PUSC of SD)}.  For each set $E$ of full measure near $x_0$ there exists a sequence
$x_j\to x_0$ with $x_j\in E$ and both $u$ and $v$ twice differentiable at $x_j$, such that }
$$
D^2_{x_j} u \ \to\ A
\and
D^2_{x_j} v \ \to\ B \qquad{\sl with } \quad A+B\leq A_0.
$$ 

\Remark{\FF.2} In particular, given an upper contact jet $(p_0, A_0)$
for the sum $w\equiv u+v$ at $x_0$, there exist two jets
$(x_0, u(x_0), p, A') \in \overline{J^+(u)}$ and 
$(x_0, v(x_0), q, B) \in \overline{J^+(v)}$ which sum to
$(x_0, w(x_0), p_0, A_0)$.  (Take $A'\equiv A+P$ where $P\equiv A_0-A-B \geq0$.)

\pf
To show that $x_0$ is an upper contact point for $u$ and $v$ we can assume that
they are both convex because of  Lemma \AA.4.  Every point is a flat lower contact point
for $v$ by the Hahn-Banach Theorem.  That is, there exists $q\in\rn$ such that
$$
v(x_0) + \bra q {y-x_0}\ \leq \ v(y) \qquad{\rm for\ all\ } y\ {\rm near\ } x_0.
\eqno{(\FF.1)}
$$
Hence, 
$$
u(y) \ =\ w(y)-v(y) \ \leq\ 
u(x_0) + \bra {p_0-q}{y-x_0} + \half \bra {A_0(y-x_0)}{y-x_0}
$$
proving that $u$ has  upper contact jet  $(p_0-q, A_0)$ at $x_0$.
The remaining first derivative results for $u$ are just restatements of Lemmas \BB.1 and \BB.2.

To prove (PUSC of SD) we apply Theorem \BB.4 to $w=u+v$ with upper contact jet $(p_0,A_0)$ at 
$x_0$ and take $E$ to be the set of points where both $u$ and $v$ are twice differentiable.
This yields a sequence $x_j\to x_0$, $x_j\in E$ with $D_{x_j}w \to p_0$
and $D^2_{x_j}w \to A_0-P$ with $P\geq0$. Moreover,
$$
D_{x_j} u\ \to\ D_{x_0}u=p,
\qquad
D_{x_j} v\ \to\ D_{x_0}v=q,
\and
p+q\ =\ p_0.
$$
In order to extract a subsequence of $\{x_j\}$ such that the second derivatives of $u$ and $v$ 
converge, i.e., 
$$
D^2_{x_j} u\ \to\ A
\and
D^2_{x_j} v\ \to\ B
$$
we prove compactness.

We can assume that $u$ and $v$ are both $\l$-quasi-convex near $x_0$ 
with the same $\l>>0$.  This implies that at each point $x\in E$ near $x_0$,
$$
-\l I   \ \leq\ D^2_x u \and -\l  I\ \leq\ D^2_x v
$$
Since $D^2_{x_j}w \to A_0-P$, we can assume that  $D^2_{x_j}w \leq A_0 - P+\e I$ for all $j$.
Hence, 
$$
D^2_{x_j}u \ =\ D^2_{x_j}w -D^2_{x_j}v \ \leq \ A_0 - P +\e I + \l I
$$
providing the upper bound for $D^2_{x_j}u$.  Hence we can assume $D^2_{x_j}u$ converges
to $A$ and similarly $D^2_{x_j}v \to B$. Moreover, $A+B = A_0+P$.\qed

\vskip .3in


\centerline{\headfont \GG.\  Strict Comparison}
 \medskip 
 
  Given a subequation $F$, recall the dual subequation $\wt F \equiv \sim(-F)$ (see [HL$_2$]).
  For quasi-convex functions, comparison always holds if one of the functions is ``strict''.
 Each subequation $G\ss \Int F$ provides an adequate  notion of $F$-strictness.

 \Theorem{\GG.1. (Strict Comparison)} {\sl
 Suppose $G$ and $F$ are subequations on a manifold $X$ with $G\ss \Int F$.
 Then  for any domain $\O\ss\ss  X$, if $u,v\in  \USC(\ob)$ are locally
 quasi-convex on $\O$ with $u\in G(\O)$ and $v\in \wt F(\O)$, comparison holds, i.e., 
  $$
 u+v\ \leq\ 0 \quad{\rm on}\ \partial \O
 \qquad\Rightarrow\qquad
 u+v\ \leq\ 0 \quad{\rm on}\  \O
 \eqno{(ZMP)}
 $$
Moreover, if $F$ and $G$ are constant coefficient on $\rn$,
then the assumption that $u$ and $v$ are quasi-convex can be dropped. }
\pf
If (ZMP) fails for $w=u+v$, then $w$ has an interior maximum point $x_0$ in $\O$,
with $m\equiv w(x_0)>0$ but $w\leq0$ on $\partial \O$.  At $x_0$, $w$ has upper contact
jet $(p_0,A_0)=(0,0)$.  Applying Theorem \FF.1 in local coordinates yields:
\medskip

(1)\ \ $(u(x_0), p, A) \in G_{x_0}$ and $(v(x_0), q, B) \in {\wt F}_{x_0}$,

\medskip

(2)\ \ $u(x_0) + v(x_0) \ =\ m  > \ 0$,

\medskip

(3)\ \ $p+q\ =\ 0$.

\medskip

(4)\ \ $A+B\ =\ -P$ with $P\geq0$.

\medskip

This contradicts $G_{x_0} \ss\Int F_{x_0}$ since 
$$
(v(x_0), q, B) \in \wt F_{x_0}  \quad\iff\quad
(-v(x_0), -q, -B) \ =\ (u(x_0)-m, p, A+P)  \notin \Int F_{x_0}  
$$
and $\Int F_{x_0}$ satisfies (N) and (P).
Again the constant coefficient part follows easily from quasi-convex approximation (cf. [HL$_1$]).
\qed

 \Note {\GG.2}  The constant coefficient case of Theorem \GG.1 was established in [HL$_2$],
 Corollary C.3,  using the Theorem on Sums ([CIL], [C]).

\vfill\eject


\centerline{\headfont   References.}
 \medskip

\item{[Al]}  A. D. Alexandrov, {\sl Almost everywhere existence of the 
second differential of a convex function and properties of 
convex surfaces connected with it (in Russian)}, 
Lenningrad State Univ. Ann. Math.  {\bf 37}   (1939),   3-35.

\smallskip

\noindent
\item{[C]}   M. G. Crandall,  {\sl  Viscosity solutions: a primer},  
pp. 1-43 in ``Viscosity Solutions and Applications''  Ed.'s Dolcetta and Lions, 
SLNM {\bf 1660}, Springer Press, New York, 1997.

 \smallskip

\item{[CIL]}   M. G. Crandall, H. Ishii and P. L. Lions {\sl
User's guide to viscosity solutions of second order partial differential equations},  
Bull. Amer. Math. Soc. (N. S.) {\bf 27} (1992), 1-67.

 \smallskip

\item{[F]}  H. Federer,   Geometric Measure Theory, Springer-Verlag, Berlin-Heidelberg, 1969.

\smallskip

\item {[HL$_{1}$]}    F. R. Harvey and H. B. Lawson, Jr.,   {\sl  Dirichlet duality and the non-linear Dirichlet problem},    Comm. on Pure and Applied Math. {\bf 62} (2009), 396-443. ArXiv:math.0710.3991

\smallskip

\item {[HL$_{2}$]}  \ \----------,   {\sl  Dirichlet duality and the nonlinear Dirichlet problem
on Riemannian manifolds},  J. Diff. Geom. {\bf 88} (2011), 395-482.   ArXiv:0912.5220.
\smallskip

\item {[HL$_{3}$]}  \ \----------,   {\sl  Characterizing the strong maximum principle},  Rendiconti di Matematica (to appear).  ArXiv:1303.1738.  ArXiv:1309.1738.
\smallskip

\item {[HL$_{4}$]}     \ \----------,   {\sl Notes on the differentiation of quasi-convex functions},     
ArXiv:1309.1772.
\smallskip

\item {[HL$_{5}$]}  \ \----------,    {\sl Tangents to subsolutions -- existence and uniqueness, Parts I and II}. 
   ArXiv:1408.5797 and  ArXiv:1408.5851.

\smallskip

\item {[HL$_{6}$]}  \ \----------,    {\sl The Dirichlet problem with prescribed interior singularities},  
Advances in Math. (to appear).     ArXiv:1508.02962.

\smallskip

\item {[S]}  Z. Slodkowski, {\sl  The Bremermann-Dirichlet problem for $q$-plurisubharmonic functions},
Ann. Scuola Norm. Sup. Pisa Cl. Sci. (4)  {\bf 11}    (1984),  303-326.

\smallskip

\item{[J]}    R. Jensen,    {\sl  The maximum principle for viscosity solutions of fully nonlinear
second order partial differential equations},   Arch. Ratl. Mech. Anal. {\bf 101}  (1988),   1-27.

\smallskip

\end